\documentclass[12pt]{amsart}
\usepackage{graphicx}
\usepackage{amsthm,mathrsfs}
\usepackage{amsfonts}
\usepackage{amsmath}
\usepackage{amssymb}

\newtheorem{theorem}{Theorem}

\newtheorem{conj}{Conjecture}
\newtheorem{definition}{Definition}
\newtheorem{example}{Example}
\newtheorem{question}{Question}
\newtheorem{problem}{Problem}

\usepackage{float}
\usepackage{url}
\usepackage{enumitem}
\usepackage{epstopdf}

\newtheorem*{remark}{Remark}

\makeatletter\def\blfootnote{\xdef\@thefnmark{}\@footnotetext}\makeatother
\allowdisplaybreaks
\title[On Pair Correlation of Sequences]{\bf On Pair Correlation of Sequences}

\author{Gerhard Larcher} 
\address{Institute of Financial Mathematics and applied Number Theory, University Linz}
\email{gerhard.larcher@jku.at}

\author{Wolfgang Stockinger} 
\address{Mathematical Institute, University of Oxford}
\email{wolfgang.stockinger@stcatz.ox.ac.uk}

\thanks{The first author is supported by the Austrian Science Fund (FWF): Project F5507-N26, which is part of the Special Research Program "Quasi-Monte Carlo Methods: Theory and Applications". The second author is supported by a special Upper Austrian grant.}

\begin{document}

\begin{abstract}
We give a survey on the concept of Poissonian pair correlation (PPC) of sequences in the unit interval, on existing and recent results and we state a list of open problems. Moreover, we present and discuss a quite recent multi-dimensional version of PPC.
\end{abstract}

\date{}
\maketitle

\section{The concept of Poissonian Pair Correlation (PPC) for sequences in $[0,1)$} \label{sect_1}
Let $x_1, x_2,\ldots,$ be a sequence of real numbers in the unit interval $[0,1)$. In the following, for some $x \in [0,1)$, we denote by $\| x \|$ the distance to the nearest integer, i.e., to be precise $\| x \|:= \min(x,1-x)$. Further, in the sequel, $\lbrace \cdot \rbrace$ will denote the fractional part of a real number. 

\begin{definition}
We say that $\left(x_n\right)_{n \geq 1} \in [0,1)$ has Poissonian pair correlation (PPC) if for all real $s > 0$ we have
$$
\underset{N \rightarrow \infty}{\lim} \frac{1}{N} \# \left\{1 \leq k \neq l \leq N \left| \left\|x_k-x_l\right\|< \frac{s}{N}\right\}\right. = 2s.
$$
\end{definition}

To put this in intuitive words, PPC means to study small distances between sequence elements, i.e., the concept of PPC deals with a ``local'' distribution property of a sequence in the unit interval.\\ \\
It is natural to expect that the pair correlation function, $R_N$, defined as 
$$
R_N(s):=\frac{1}{N} \# \left\{1 \leq k \neq l \leq N \left|\left\|x_k-x_l\right\|< \frac{s}{N} \right\}\right.
$$
tends to $2s$. We give the following heuristic explanation for this limit behaviour:

\begin{center}
\includegraphics[angle=0,width=100mm]{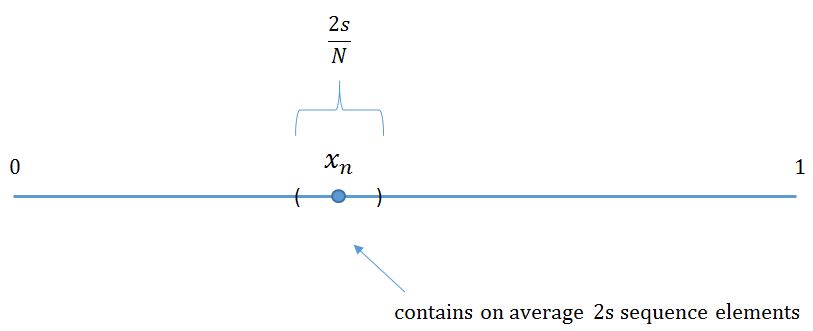}
~\\
Figure 1
\end{center}
~\\
Consider a fixed $N$, and fix a sequence element $x_n$ for some $1\leq n\leq N$. Then, the region around $x_n$ with length $\frac{2s}{N}$ (see Figure 1) i expected to contain $2s\frac{N-1}{N}$ of the remaining $(N-1)$ points $x_i$, for $i=1,\ldots, N$ and $i \neq n$.\\
Consequently, this means, on average there are $2s \frac{N-1}{N}$ different indices $i = 1,\ldots, N$ with $i \neq n$, such that 
\begin{equation*}
\left\|x_i - x_n\right\| < \frac{s}{N}.
\end{equation*}
Since $n$ can attain values between $1$ and $N$, we expect that there are $2s(N-1)$ pairs with 
\begin{equation*}
\left\|x_k-x_l\right\|< \frac{s}{N}, \qquad \text{for } 1 \leq k  \neq l \leq N.
\end{equation*}
Hence, we expect the quantity $R_N(s)$ to be approximately $2s \frac{N-1}{N}$, and therefore
\begin{equation*}
\underset{N \rightarrow \infty}{\lim} R_N (s) = 2s.
\end{equation*}
Indeed, it can be shown that, in a certain sense, almost every sequence $x_1, x_2, \ldots$ in $[0,1)$ has PPC. To be precise, if we consider a sequence $(X_n)_{n \geq 1}$ of i.i.d.~ random variables drawn from the uniform distribution on $[0,1)$, then 
\begin{equation*}
\underset{N \rightarrow \infty}{\lim} \frac{1}{N} \# \left\{1 \leq k \neq l \leq N \left| \left\|X_k-X_l\right\|< \frac{s}{N}\right\}\right. = 2s, 
\end{equation*}
almost surely. \\ \\
We want to briefly mention that the original motivation for the investigation of the PPC property comes from quantum physics. Roughly speaking, the concept is related to the distribution properties of the discrete energy spectrum $\lambda_1, \lambda_2, \ldots$ of a Hamiltonian operator of a quantum system. The famous Berry-Tabor-Conjecture in quantum physics now states, that this discrete energy spectrum (ingoring degenerate cases) has PPC. For more details on the connection to quantum physics, we refer the reader to, for example, \cite{AAL} and the references cited therein.\\ \\
Note that for several quantum systems the discrete energy spectrum $\lambda_1,\lambda_2,\ldots$ has the following special form
$$
\left(\lambda_n\right)_{n \geq 1} = \left(\left\{a_n \alpha\right\}\right)_{n \geq 1},
$$
where $\alpha$ is a real constant, and $\left(a_n\right)_{n \geq 1}$ is a given sequence of positive integers. Therefore, in the 1990's Rudnick, Sarnak and Zaharescu started to investigate the PPC property of sequences of the form $\left(\left\{a_n \alpha\right\}\right)_{n \geq 1}$ in $[0,1)$ from a purely mathematical point of view. \\
Whenever, in the following, we consider such sequences we restrict the setting to strictly increasing sequences of positive integers. The most basic example of such a sequence is the classical Kronecker sequence $\left(\left\{n \alpha\right\}\right)_{n \geq 1}$. This sequence does not have the PPC property for any choice of $\alpha$. In most of the seminal papers on PPC, this fact was argued by taking the famous Three-Gap-Theorem into account (see e.g., \cite{not1, not2, not3}). \\
The Three-Gap-Theorem states the following: For every choice of $\alpha$ and for every $N$, the gaps between neighbouring points of the set
\begin{equation*}
\left\{1 \alpha\right\}, \left\{2 \alpha\right\}, \ldots, \left\{N \alpha\right\}
\end{equation*}
can have at most three different lengths. A sequence with such a gap structure does not exhibit a random behaviour and therefore it is reasonable to expect that it cannot have PPC. Nonetheless, it is not immediately clear that this argument is indeed valid. \\

\begin{center}
\includegraphics[angle=0,width=100mm]{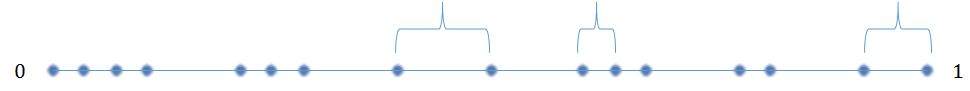}
~\\
Figure 2
\end{center}
~\\
To argue that, we want to emphasize that the elements of a sequence satisfying such a weak gap structure could be ordered in a way, such that ``many different'' distances between (not necessarily neighbouring) elements can occur (see Figure 2). However, for the Kronecker sequence a very simple argument can be given to deduce the fact that it does not have PPC for any choice of $\alpha$:  \\ \\
Let $\alpha \sim \frac{p_n}{q_n}$ where $\frac{p_n}{q_n}$ is a best approximation fraction to $\alpha$ with $\left(p_n,q_n\right)=1$. It is well-known from basic Diophantine approximation theory that $\alpha = \frac{p_n}{q_n} + \theta_n$, with either $0\leq \theta_n < \frac{1}{2q_n^2}$, or $-\frac{1}{2q_n^2} < \theta < 0$.\\
Let us assume the first case. Then, the set of points 
\begin{equation*}
\left\{1 \alpha\right\}, \left\{2 \alpha\right\}, \ldots, \left\{N \alpha\right\}
\end{equation*}
equals the set of points 
\begin{equation*}
\frac{0}{N} + \varphi_0, \frac{1}{N}+\varphi_1, \ldots, \frac{N-1}{N} + \varphi_{N-1},
\end{equation*}
with $0 \leq \varphi_i < \frac{1}{2N}$, for $i=0,\ldots,N-1$ (see red points in Figure 3).

\begin{center}
\includegraphics[angle=0,width=100mm]{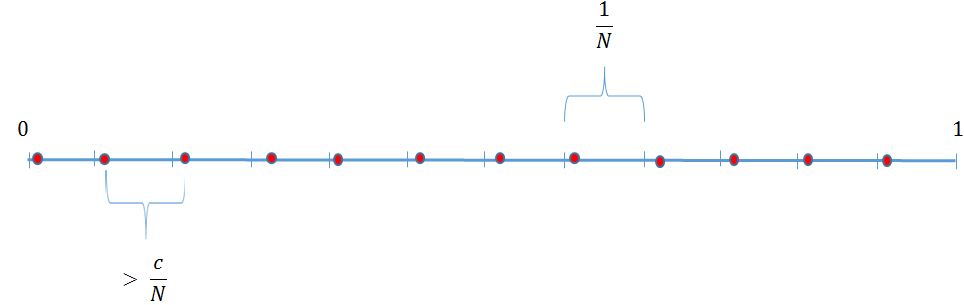}
~\\
Figure 3
\end{center}
~\\
Hence, two arbitrary elements of this point set have a distance of at least $\frac{1}{2N}$. Thus, for the choice $s=\frac{1}{4}$, we get that $R_N (s)=0$ and consequently the pair correlation function cannot tend to $2s =\frac{1}{2}$ for $N$ to infinity.\\ \\
In fact, a deeper investigation reveals that the Three-Gap-Theorem is indeed a valid argument to deduce the result on the PPC strucuture of the Kronecker sequence. It is an immediate consequence of the following theorem, which was proven in \cite{not4}, in combination with the Three-Gap-Theorem.

\begin{theorem}\label{th:gap}
Let $\left(x_n\right)_{n \geq 1}$ be a ``weak finite-gap-sequence'', i.e., there exists an integer $L$ and indices $N_1<N_2<N_3< \ldots$ such that for all $i$ the set $x_1, x_2, \ldots, x_{N_{i}}$ has at most $L$ different gap lengths between neighbouring elements. Then, $\left(x_n\right)_{n \geq 1}$ does not have PPC.
\end{theorem}

Let us now come to ``posititve results'' and to the study of the metrical pair correlation theory of sequences of the form $(\lbrace a_n \alpha \rbrace)_{n \geq 1}$. In \cite{not3} Rudnick and Sarnak showed the following:

\begin{theorem}
The sequence $\left(\left\{n^d\alpha\right\}\right)_{n \geq 1}$ with an integer $d \geq 2$ has PPC for almost all $\alpha$.
\end{theorem}
The case $d=2$ is of particular interest, as in this setting the spacings of sequence elements are related to the distances between the energy levels of the so-called  ``boxed oscillator'', i.e., the study of the PPC property of $\left(\left\{n^2\alpha\right\}\right)_{n \geq 1}$ is of special importance in quantum physics. The PPC property and also the gap distribution of the sequence $\left(\left\{n^2\alpha\right\}\right)_{n \geq 1}$ was further investigated by several authors (see \cite{not5,not6,not2}) and they could also derive the metrical result for $d=2$. Heath-Brown could even show slightly more:  

\begin{theorem}
The sequence $\left(\left\{n^2\alpha\right\}\right)_{n \geq 1}$ has PPC for almost all real numbers $\alpha$. Moreover, there is a dense set of constructible values of $\alpha$ for which the PPC property holds, i.e., there is an informal algorithm, which, for any closed interval $I$ of positive length, provides a convergent sequence of rational numbers belonging to $I$, whose limit $\alpha$ satisfies that $\left(\left\{n^2\alpha\right\}\right)_{n \geq 1}$ has PPC.
\end{theorem}

These three results are only metrical statements. However, until now, however, no single explicit $\alpha$ is known such that $\left(\left\{n^2 \alpha\right\}\right)_{n \geq 1}$ (or $\left(\left\{n^d \alpha\right\}\right)_{n \geq 1}$ for any integer $d \geq 2$) has PPC. Nonetheless, we know that it is \textbf{not true}, that $\left(\left\{n^2 \alpha\right\}\right)_{n \geq 1}$ has PPC for \textbf{all} irrational $\alpha$. Consider the following example of an $\alpha$ that is in a certain sense well-approximable: If $\alpha$ is an irrational number such that $\left|\alpha - \frac{a}{q}\right| < \frac{1}{4 q^3}$ for infinitely many integers $a$ and $q$, then $\left(\left\{n^2\alpha\right\}\right)_{n \geq 1}$ does not have PPC (see \cite{not5}).  \\ \\
On the other hand, it is conjectured that for an $\alpha$ which is not too well approximable, we have PPC for $\left(\left\{n^2\alpha\right\}\right)_{n \geq 1}$. To be precise: Let $\alpha$ be such that for every $\varepsilon > 0$ there is a $c(\varepsilon)>0$ with $\left|\alpha - \frac{a}{q}\right|> c(\varepsilon) \frac{1}{q^{2+\varepsilon}}$ for all $a,q \in \mathbb{Z}$, then $\left(\left\{n^2 \alpha\right\}\right)_{n \geq 1}$ has PPC (see e.g., \cite{not5}). This property for an irrational $\alpha$ is often referred to as Diophantine. It is well-known that almost all irrationals are Diophantine, e.g., every real irrational algebraic number has this property. Above discussion illustrates that the pair correlation theory of sequences $\left(\left\{n^d\alpha\right\}\right)_{n \geq 1}$ is strongly related to the Diophantine properties of $\alpha$.  \\ \\
The case of lacunary sequences $\left(a_n\right)_{n \geq 1}$ was considered for example by Rudnick and Zaharescu (see \cite{not1}) or by Berkes, Philipp and Tichy (see \cite{not7}). We recall that a sequence $\left(a_n\right)_{n \geq 1}$ is a lacunary sequence if there exists a $c > 1$ such that $\frac{a_{n+1}}{a_n} > c$ for all $n \geq N(c)$. Again, they obtained the following metrical result:

\begin{theorem}
Let $\left(a_n\right)_{n \geq 1}$ be a lacunary sequence of positive integers. Then $\left(\left\{a_n \alpha\right\}\right)_{n \geq 1}$ has PPC for almost all $\alpha$.
\end{theorem}

We may again ask for explicit examples of lacunary  sequences $\left(a_n\right)_{n \geq 1}$ and $\alpha \in \mathbb{R}$ such that $\left(\left\{a_n \alpha\right\}\right)_{n \geq 1}$ has PPC.\\
One of the most basic examples of a lacunary sequence of integers certainly is the sequence $\left(2^n\right)_{n \geq 1}$. In a first step, we may restrict the possible candidates $\alpha$ for which $\left(\left\{2^n \alpha\right\}\right)_{n \geq 1}$ could have PPC. To do so, we consider the following result which has been shown independently by Grepstad and Larcher \cite{not8}, Aistleitner, Lachmann and Pausinger \cite{not9}, and Steinerberger \cite{not10}:

\begin{theorem}
If the sequence $\left(x_n\right)_{n \geq 1}$ in $[0,1)$ has PPC, then $\left(x_n\right)_{n \geq 1}$ is uniformly distributed in $[0,1)$.
\end{theorem}

\begin{remark}
The paper of Grepstad and Larcher also contains a quantitative version of this result. Roughly speaking: If $R_N(s)$ tends to $2s$ ``fast in some uniform sense'', then the discrepancy $D_N$ of the sequence $\left(x_n\right)_{n \geq 1}$ cannot tend to zero ``too slowly''.
\end{remark}
Having this result in mind, we can restrict the set of possible choices for $\alpha$, such that $\left(\left\{2^n \alpha\right\}\right)_{n \geq 1}$ has PPC. The above theorem implies that for such an $\alpha$ the sequence $\left(\left\{2^n \alpha\right\}\right)_{n \geq 1}$ has to be uniformly distributed, and hence $\alpha$ has to be normal in base 2.\\
The most well-known example of a real $\alpha$ which is normal in base 2 is the Champernowne number in base 2, i.e., the number $\alpha$ which has in base 2 the digit representation
$$
\alpha = 0.~01~10~11~100~101~110~111~1000\ldots
$$
However, it was shown by Pirsic and Stockinger in \cite{not11}, that for $\alpha$ the Champernowne number, the sequence $\left(\left\{2^n \alpha\right\}\right)_{n \geq 1}$ does not have PPC. Also for further concrete examples like Stoneham-numbers or infinite de Bruijn-words, the sequence $\left(\left\{2^n\alpha\right\}\right)_{n \geq 1}$ does not have PPC (see \cite{not4}).\\ \\
Indeed, until now we do not know any concrete example of $\left(a_n\right)_{n \geq1}$, a lacunary sequence, and a real $\alpha$ such that $\left(\left\{a_n \alpha\right\}\right)_{n \geq 1}$ does have PPC.\\ \\
Recently, a much more general metric result on PPC of sequences of the form $\left(\left\{a_n \alpha\right\}\right)_{n \geq 1}$ was given in \cite{not12} which shows that there is an intimate connection between the concept of PPC of sequences $\left(\left\{a_n \alpha\right\}\right)_{n \geq 1}$ and the notion of additive energy of the sequence $\left(a_n\right)_{n \geq 1}$. The concept of additive energy plays a central role in additive combinatorics and also appears in the study of the metrical discrepancy theory of sequences $\left(\left\{a_n \alpha\right\}\right)_{n \geq 1}$ (see \cite{not26,not25}).\\
For a strictly increasing sequence $a_1<a_2 < a_3 < \ldots$ of positive integers, we consider the first $N$ elements $a_1, \ldots, a_N$. The additive energy of $a_1, \ldots, a_N$ is given by
$$
E\left(a_1, \ldots,a_N\right) := \sum_{\underset{a_i-a_j=a_k-a_l}{1 \leq i, j,k,l \leq N}} 1.
$$
It is obvious that $N^2 \leq E\left(a_1,\ldots,a_N\right) \leq N^3$ always holds. In \cite{not12} the following was shown:

\begin{theorem}
Let $\left(a_n\right)_{n \geq1}$ be a strictly increasing sequence of integers such that there exists $\varepsilon > 0$ with
$$
E\left(a_1, \ldots, a_N\right) = \mathcal{O} \left(N^{3-\varepsilon}\right),
$$
then $\left(\left\{a_n \alpha\right\}\right)_{n \geq 1}$ has PPC for almost all $\alpha$.
\end{theorem}
This result recovers all above mentioned metrical results and implies several new results and examples.

\begin{example}
If $\left(a_n\right)_{n \geq 1}$ is lacunary, then $E\left(a_1, \ldots,a_N\right) = \mathcal{O} \left(N^2\right)$, hence $\left(\left\{a_n \alpha\right\}\right)_{n \geq 1}$ has PPC for almost all $\alpha$.
\end{example}

\begin{example}
If $\left(a_n\right)_{n \geq 1}$ are the values of a polynomial $f(n) \in \mathbb{Z} \left[x\right]$ of degree $d \geq 2$, then $E\left(a_1, \ldots, a_N\right) = \mathcal{O} \left(N^{2 + \varepsilon}\right)$ for all $\varepsilon > 0$, hence $\left(\left\{a_n \alpha\right\}\right)_{n \geq 1}$ has PPC for almost all $\alpha$.
\end{example}

\begin{example}
Let $\left(a_n\right)_{n \geq 1}$ be a convex sequence, i.e., $a_n -a_{n-1} < a_{n+1}-a_n$ for all $n$, then it was shown by Konjagin \cite{not13}, that $E\left(a_1, \ldots,a_N\right) = \mathcal{O} \left(N^{\frac{5}{2}}\right)$, hence $\left(\left\{a_n \alpha\right\}\right)_{n \geq 1}$ has PPC for almost all $\alpha$.
\end{example}

\begin{example}
If $a_n = \left[\beta  n^c\right]$ for some $\beta > 0$ and $c > 1$, then 
$$
E\left(a_1,\ldots,a_N\right) = \mathcal{O} \left(\max\left(N^{\frac{5}{2}}, N^{4-c}\right)\right),
$$
hence $\left(\left\{a_n \alpha\right\}\right)_{n \geq 1}$ has PPC for almost all $\alpha$ (see \cite{not14}).
\end{example}
The above theorem immediately raises two natural questions:

\begin{question} \label{qu_a}
Is it possible for an increasing sequence of distinct integers $(a_n)_{n \geq 1}$ which satisfies $E(a_1, \ldots, a_N) = \Omega(N^3)$ that the sequence $(\lbrace a_n \alpha \rbrace)_{n \geq 1}$ has PPC for almost all $\alpha$?
\end{question}

\begin{question} \label{qu_b}
If, for almost all $\alpha$, $\left(\left\{a_n \alpha\right\}\right)_{n \geq 1}$ does \textbf{not} have PPC, does this imply $E\left(a_1,\ldots,a_N\right) = \Omega \left(N^3\right)$?
\end{question}

Both questions were answered by J. Bourgain in an Appendix to \cite{not12}.\\ \\
Concerning Question~\ref{qu_a}, Bourgain showed:
\begin{itemize}
\item If $E\left(a_1, \ldots,a_N\right) = \Omega \left(N^3\right)$, then there exists a set of positive measure such that $\left(\left\{a_n \alpha\right\}\right)_{n \geq 1}$ does not have PPC for every $\alpha$ in this set.
\end{itemize}
This result was improved by Lachmann and Technau \cite{not15}:
\begin{itemize}
\item If $E\left(a_1, \ldots,a_N\right) = \Omega \left(N^3\right)$, then there exists a set of full Lebesgue measure such that $\left(\left\{a_n \alpha\right\}\right)_{n \geq 1}$ does not have PPC for every $\alpha$ contained in this set.
\end{itemize}
Finally, in \cite{not16} this result was improved to its final form:

\begin{theorem}
If $E\left(a_1, \ldots, a_N\right) = \Omega \left(N^3\right)$, then there is \textbf{no} $\alpha$ such that $\left(\left\{a_n \alpha\right\}\right)_{n \geq 1}$ has PPC.
\end{theorem}

Concerning Question~\ref{qu_b} Bourgain showed that the answer to this question is \textbf{``no''}: He gave a construction for  a sequence $\left(a_n\right)_{n\geq 1}$ with $E\left(a_1, \ldots,a_N\right) = o\left(N^3\right)$, such that $\left(\left\{a_n \alpha\right\}\right)_{n \geq 1}$ does not have PPC for almost all $\alpha$.\\ \\
Up to now, we have the following situation:\\
$E\left(a_1, \ldots,a_N\right) = \Omega \left(N^3\right)$ implies that there is \textbf{no} $\alpha$ such that the sequence $\left(\left\{a_n \alpha\right\}\right)_{n \geq 1}$ has PPC. \\
$E\left(a_1, \ldots,a_N\right) = \mathcal{O}\left(N^{3-\varepsilon}\right)$ for some $\varepsilon > 0$ implies PPC for almost all $\alpha$. \\ \\
The result by Aistleitner, Larcher and Lewko, was first extended by Bloom, Chow, Gafni and Walker, albeit under an additional density condition on the integer sequence $(a_n)_{n \geq 1}$ (see \cite{not27}). 
\begin{theorem}
Let $a_1, \ldots, a_N$ be the first $N$ elements of an increasing sequence of positive integers $( a_n )_{n \geq 1}$ satisfying the following density condition:
\begin{equation*}
\delta(N) = \Omega_{\varepsilon} \left( \frac{1}{(\log N)^{2 +2\varepsilon}} \right), 
\end{equation*}
where $\delta(N):=N^{-1} \#( \lbrace a_1, \ldots, a_N \rbrace \cap \lbrace 1, \ldots, N \rbrace)$ and suppose that 
\begin{equation*}
E(a_1, \ldots, a_n) = \mathcal{O}_{\varepsilon}\left( \frac{N^3}{(\log N)^{2+ \varepsilon} } \right),
\end{equation*}
for some $\varepsilon>0$, then, for almost all $\alpha$, the sequence $(\lbrace a_n \alpha \rbrace)_{n \geq 1}$ has PPC.
\end{theorem}
Recently, Bloom and Walker (see \cite{not17}) improved over this result by showing the following theorem. 
\begin{theorem}\label{TH:THWalk}
There exists an absolute positive constant $C$ such that the following is true. Suppose that
\begin{equation*}
E(a_1, \ldots, a_N) = \mathcal{O}\left(\frac{N^3}{(\log N)^C} \right),
\end{equation*} 
then for almost all $\alpha$, $(\lbrace a_n \alpha \rbrace)_{n \geq 1}$ has PPC.
\end{theorem}
A consequence of this result is the following theorem:
\begin{theorem}
Let $(a_n)_{n \geq 1}$ be an arbitrary infinite subset of the squares. Then $(a_n)_{n \geq 1}$ is metric Poissonian, i.e., for almost all $\alpha$, $(\lbrace a_n \alpha \rbrace)_{n \geq 1}$ has PPC. 
\end{theorem}
To see that this result is valid, we note that, if $a_1, \ldots, a_N$ denotes a finite set of squares, then $E(a_1, \ldots, a_N)= \mathcal{O}( N^3 \exp({-c_1 \log^{c_2} N)})$ for some absolute positive constants $c_1$ and $c_2$, see e.g., \cite{not24}. \\ \\
The proof of Theorem \ref{TH:THWalk} relies on a new bound for GCD sums with $\alpha=1/2$, 
which improves over the bound by Bondarenko and Seip (see \cite{not18}), 
if the additive energy of $a_1, \ldots, a_N$ is sufficiently large.
Note that the constant $C$ was not specified in the above mentioned paper, 
but the authors thereof conjecture that Theorem~\ref{TH:THWalk} holds for $C>1$ already. This result would be best possible. To see this, consider the following result by Walker \cite{not19}:

\begin{theorem}
Let $\left(a_n\right)_{n \geq 1} = \left(p_n\right)_{n \geq 1}$ be the sequence of primes (note that for the primes we have $E\left(p_1, \ldots,p_N\right) \asymp \frac{N^3}{\log  N}$). Then, $\left(\left\{p_n \alpha \right\}\right)_{n \geq 1}$ does not have PPC for almost all $\alpha$.
\end{theorem}

The region between $\mathcal{O}\left(\frac{N^3}{(\log N)^C} \right)$, $C>1$, and $\Omega\left(N^3\right)$ is therefore the interesting region and one might speculate about a sharp threshold which allows to fully describe the metrical pair correlation theory in terms of the additive energy. Further constructions and examples of sequences in this ``interesting region'', with an even smaller additive energy compared to the primes, were given by Lachmann and Technau \cite{not15}:

\begin{theorem}
There exists a strictly increasing sequence of positive integers $\left(a_n\right)_{n \geq1}$ with 
$$
E\left(a_1, \ldots,a_N\right)=\mathcal{O} \left(\frac{N^3}{\log N  (\log \log N)}\right)
$$
such that $\left(\left\{a_n \alpha\right\}\right)_{n \geq 1}$ does not have PPC for almost all $\alpha$.
\end{theorem}

On the other hand they gave a positive result of the following form:

\begin{theorem}
There exists a strictly increasing sequence of positive integers $\left(a_n\right)_{n \geq 1}$ with 
$$
E\left(a_1, \ldots,a_N\right) = \Omega \left(\frac{N^3}{\log N  \left(\log \log\right)^{1+\varepsilon}}\right)
$$ 
for all $\varepsilon > 0$, such that $\left(\left\{a_n \alpha\right\}\right)_{n \geq 1}$ has PPC for almost all $\alpha$.
\end{theorem}

We have the following situation as illustrated in Figure 4.

\begin{center}
\includegraphics[angle=0,width=1\textwidth]{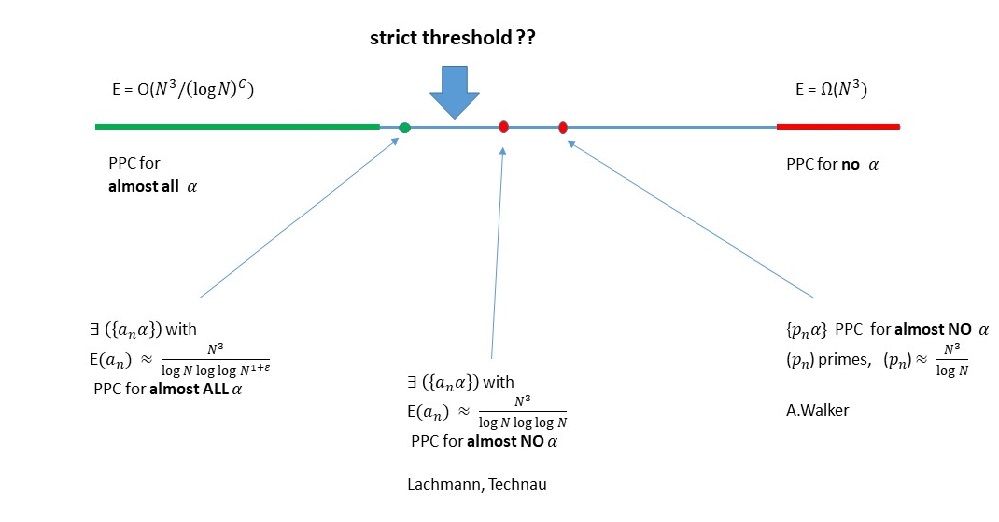}
~\\
Figure 4
\end{center}
~\\
The following question therefore is near at hand:\\ \\
Is there a strict threshold $T$ such that an additive energy of magnitude smaller than $T$ implies PPC of $\left(\left\{a_n \alpha\right\}\right)_{n \geq 1}$ for almost all $\alpha$ and an additive energy of magnitude larger than $T$ implies PPC for $\left(\left\{a_n \alpha\right\}\right)_{n \geq 1}$ for almost no $\alpha$? The fundamental question concerning such a putative threshold was raised in \cite{not17}. The authors of this paper conjectured that there is a sharp Khintchine-type threshold, i.e., if $E (a_1, \ldots, a_N) =\Theta( N^3\psi(N))$, for some weakly decreasing function $\psi : \mathbb{Z}_{\geq 1} \to [0,1]$, then, for almost all $\alpha$, $(\lbrace a_n \alpha \rbrace)_{n \geq 1}$ has PPC if and only if
\begin{equation*}
\sum_{N \geq 1} \frac{\psi(N)}{N}
\end{equation*}
converges. \\ \\
The negative answer to this question was given by Aistleitner, Lachmann and Technau \cite{not20}:

\begin{theorem}
There exists a sequence $\left(a_n\right)_{n \geq 1}$ of integers with 
$$
E\left(a_1, \ldots, a_N\right) = \Omega \left(\frac{N^3}{\left(\log N\right)^{\frac{3}{4}+\varepsilon}}\right)
$$ 
such that $\left(\left\{a_n \alpha\right\}\right)_{n \geq 1}$ has PPC for almost all $\alpha$. Hence a threshold $T$ cannot exist.
\end{theorem}
To conclude, the additive energy is not enough to fully describe the metrical pair correlation theory. Some further number theoretic properties need to be considered to cope with that problem.
\section{The concept of Poissonian Pair Correlation (PPC) for sequences in $[0,1)^d$}
Of course it makes sense to generalize the concept of PPC to the multi-dimensional setting. One way to generalize the one-dimensional concept to a multi-dimensioanl setting was defined and discussed in \cite{not21} (for a more general analysis of a multi-dimensional PPC concept, we refer to the recent work \cite{not29}). Here, we present the definition of \cite{not21}. 
\begin{definition}
Let $\left(x_n\right)_{n \geq 1}$ be a sequence in the $d$-dimensional unit-cube $\left[\left.0,1\right.\right)^d$. We say that $\left(x_n\right)_{n \geq 1}$ has PPC if for all $s > 0$ we have
$$
\underset{N \rightarrow \infty}{\lim} \frac{1}{N} \# \left\{1 \leq k \neq l \leq N \left| \left\|x_k-x_l\right\|_{\infty} < \frac{s}{N^{\frac{1}{d}}}\right\}\right. = \left(2s\right)^d.
$$
\end{definition} 
For this definition of $d$-dimensional PPC it again follows that $\left(x_n\right)_{n \geq 1}$ with PPC is uniformly distributed in $\left[\left.0,1\right.\right)^d$. Moreover, for many of the above mentioned results in dimension $d=1$ we have analogous statements in dimension $d \geq 2$. For example the $d$-dimensional Kronecker sequence 
\begin{equation*}
\left(\left\{n \alpha_1\right\}, \left\{n \alpha_2\right\}, \ldots,\left\{n \alpha_d\right\}\right)_{n \geq 1}
\end{equation*}
never has PPC. The proof of this fact however needs a bit more subtle arguments than in dimension $1$. \\ \\
Naturally, we would also expect that under the same condition on the additive energy as in Theorem \ref{TH:THWalk}, the sequence 
\begin{equation*}
( \lbrace a_n \boldsymbol{\alpha} \rbrace)_{n \geq 1}
\end{equation*}
has Poissonian pair correlations for almost all instances and, in fact, we have the following even better result, which is a consequence of better bounds on GCD sums for larger exponents than $1/2$:
\begin{theorem}\label{TH:THMet}
Let $a_1, \ldots, a_N$ denote the first $N$ elements of $(a_n)_{n \geq 1}$ and suppose that 
\begin{equation*}
E(a_1, \ldots, a_N) = \mathcal{O}\left(\frac{N^3}{(\log N)^{1+ \varepsilon}} \right), \text{ for any } \varepsilon > 0,
\end{equation*} 
then for almost all choices of $\boldsymbol{\alpha}=(\alpha_1, \ldots, \alpha_d) \in \mathbb{R}^d$,
\begin{equation*}
( \lbrace a_n \boldsymbol{\alpha} \rbrace)_{n \geq 1}
\end{equation*}
has PPC.
\end{theorem}
However, if the additive energy is of maximal order, i.e., if we have $E(a_1, \ldots, a_N)=\Omega(N^3)$, then there is no $\boldsymbol{\alpha}$ such that $(\lbrace a_n \boldsymbol{\alpha} \rbrace)_{n \geq 1}$ has PPC:
\begin{theorem}\label{TH:THMAX}
If  $E(a_1, \ldots, a_N)=\Omega(N^3)$, then for any choice of $\boldsymbol{\alpha} =(\alpha_1, \ldots, \alpha_d) \in \mathbb{R}^d$ the sequence
\begin{equation*}
( \lbrace a_n \boldsymbol{\alpha} \rbrace)_{n \geq 1},
\end{equation*} 
does not have Poissonian pair correlations. 
\end{theorem}
\section{Open Problems}
Many questions related to the concept of PPC are still open, and we will state some of them in this section as open problems.
\begin{problem}
Is it possible to extend the ``green region'' and/or the ``red region'' in Figure~4? That means: Are there functions $\varphi(n)$ (which increases slower than $(\log N)^C$, for $C$ the constant in Theorem \ref{TH:THWalk}) and $\psi (n)$ both tending to $+\infty$ for $n$ to infinity, such that:\\
If $E\left(a_1, \ldots,a_N\right) = \mathcal{O} \left(\frac{N^3}{\varphi(N)}\right)$, then, for almost \textbf{all} $\alpha$, $\left(\left\{a_n \alpha\right\}\right)_{n \geq 1}$ has PPC.\\
If $E\left(a_1, \ldots,a_N\right) = \Omega \left(\frac{N^3}{\psi(N)}\right)$, then there is \textbf{no} $\alpha$ such that $\left(\left\{a_n \alpha\right\}\right)_{n \geq 1}$ has PPC.
\end{problem}

\begin{problem}
We know that if $E\left(a_1, \ldots,a_N\right) = \Omega \left(N^3\right)$, then $\left(\left\{a_n \alpha\right\}\right)_{n \geq 1}$ has PPC for \textbf{no} $\alpha$. We consider the following question to be of high interest: Is there a sequence $\left(a_n\right)_{n \geq 1}$ with the property that for almost all $\alpha$, $\left(\left\{a_n \alpha\right\}\right)_{n \geq 1}$ does not have PPC, but, there exists a set of zero measure such that $\left(\left\{a_n \alpha\right\}\right)_{n \geq 1}$ has PPC for every $\alpha$ contained in this set?
\end{problem}

Indeed, we believe, that this is not possible, i.e.,:

\begin{conj}
If, for almost all $\alpha$, $\left(\left\{a_n \alpha\right\}\right)_{n \geq1}$ does not have PPC, then it has PPC for \textbf{no} $\alpha$.
\end{conj}

For example (by the result of A. Walker) this would imply: $\left(\left\{p_n \alpha\right\}\right)_{n \geq 1}$ has PPC for \textbf{no} $\alpha$.

\begin{problem}
Although the metrical theory of sequences of the form $( \lbrace a_n \boldsymbol{\alpha} \rbrace)_{n \geq 1}$ seems to be well-established, we do not know any explicit construction of $\boldsymbol{\alpha}$ (not even in the one-dimensional case) such that $( \lbrace a_n \boldsymbol{\alpha} \rbrace)_{n \geq 1}$ has Poissonian pair correlations. It is in general very hard to construct sequences on the torus having the PPC property. The only known explicit examples -- to the best of our knowledge -- of sequences with this property are $\lbrace \sqrt{n} \rbrace_{n \geq 1}$ (see \cite{not22}) and certain directions of vectors in an affine Euclidean lattice (see \cite{not23}). Hence, of course, it would be of high interest to find more concrete examples of sequences with PPC.
\end{problem}

\begin{problem}
This problem concerns a possible extension of Theorem \ref{th:gap} mentioned above. We recall, that we have shown that a sequence $\left(x_n\right)_{n \geq 1}$ with a weak finite gap property never has PPC. We wonder whether this result can be improved by showing that it still holds if we have to deal with a sequence with not necessarily ``finite gap-'' but with a ``slowly-growing-gap-'' property, i.e., with a sequence with the following property:\\
There is a (slowly growing) function $L$ and a sequence of indices $N_1<N_2<N_3<\ldots$ such that $x_1,x_2,\ldots,x_{N_{i}}$ always has gaps of at most $L\left(N_i\right)$ different lengths.
\end{problem}

\begin{problem}
Let $\left(\boldsymbol{x}_n\right)_{n \geq 1}$ be the $d$-dimensional Halton-sequence in any bases $q_1, \ldots,q_d$, where $d \geq 2$. Does $\left(\boldsymbol{x}_n\right)_{n \geq 1}$ have PPC or not? Of course, we strongly conjecture that it does not have PPC. In dimension $d=1$ the Halton-sequence is the well-known van der Corput sequence. In this case the PPC property trivially does not hold. In fact, it should not be too hard to prove this in the multidimensional case, too.
\end{problem}

The last problem we want to state concerns a multi-dimensional version of the metrical PPC result for the primes.

\begin{problem}
Is it true that for almost all instances of $\boldsymbol{\alpha}$ the sequence $( \lbrace p_n \boldsymbol{\alpha} \rbrace)_{n \geq 1}$, where $(p_n)_{n \geq 1}$ denotes the primes, does not have PPC? 
\end{problem}

\end{document}